\documentclass[12pt, a4paper, twoside]{article}
\usepackage{amssymb}
\usepackage[leqno]{amsmath}
\usepackage{latexsym}
\usepackage{enumerate}
\usepackage{amsmath, amssymb, amscd}
\Roman{equation}

\DeclareMathOperator{\rk}{rk}
\DeclareMathOperator{\Pic}{Pic}

\begin{document}

\title{\textbf{\Large{ACM bundles of rank 2 on quartic hypersurfaces in $\mathbb{P}^3$ and Lazarsfeld-Mukai bundles}}}

\author{Kenta Watanabe \thanks{Nihon University, College of Science and Technology,   7-24-1 Narashinodai Funabashi city Chiba 274-8501 Japan , {\it E-mail address:watanabe.kenta@nihon-u.ac.jp}, Telephone numbers: 090-9777-1974} }

\date{}

\maketitle 


\begin{abstract}

\noindent Let $X$ be a smooth quartic hypersurface in $\mathbb{P}^3$. By the Brill-Noether theory of curves on K3 surfaces, if a rank 2 aCM bundle on $X$ is globally generated, then it is the Lazarsfeld-Mukai bundle $E_{C,Z}$ associated with a smooth curve $C$ on $X$ and a base point free pencil $Z$ on $C$. In this paper, we will focus on the classification of such bundles on $X$ to investigate aCM bundles of rank 2 on $X$. Concretely, we will give a necessary condition for a rank 2 vector bundle of type $E_{C,Z}$ to be indecomposable initialized and aCM, in the case where the class of $C$ in $\Pic(X)$ is contained in the sublattice of rank 2 generated by the hyperplane class of $X$ and a non-trivial initialized aCM line bundle on $X$.

\end{abstract}

\section{Introduction} It is interesting to consider the problem concerning the splitting of vector bundles on smooth hypersurfaces in $\mathbb{P}^n$ in the point of view of the Hartshorne conjecture. Let $X$ be a smooth hypersurface in $\mathbb{P}^n$, and let $H$ be a hyperplane section of it. Then we call a vector bundle $E$ on $X$ an {\it{arithmetically Cohen-Macaulay}} ({{\it{aCM}} for short) bundle if $H^i(E\otimes \mathcal{O}_X(lH))=0$ for all integers $l\in\mathbb{Z}$ and $1\leq i\leq n-2$. If $X=\mathbb{P}^{n-1}$, any indecomposable aCM bundle on $X$ is a line bundle ([12]). If $X$ is a smooth quadric, then any non-split aCM bundle of rank $r\geq2$ on $X$ is isomorphic to a spinor bundle up to twisting by the hyperplane class of $X$ ([1]). By these facts, any aCM bundle on a smooth hypersurface $X$ of degree $d\leq2$ in $\mathbb{P}^3$ splits into line bundles. However, we can find indecomposable aCM bundles on smooth hypersurfaces of degree $d\geq3$ in $\mathbb{P}^3$. If a vector bundle $E$ on $X$ is aCM, then $E\otimes\mathcal{O}_X(-H)$ is also aCM. Hence, previously, an aCM bundle $E$ on $X$ of degree $d\geq3$ in $\mathbb{P}^3$  satisfying $h^0(E)\neq0$ and $h^0(E\otimes\mathcal{O}_X(-H))=0$ which is called an {\it{initialized}} aCM bundle is often studied by many people. For example, in the case where $d=3$, Casanellas and Hartshorne ([3]) have constructed families of indecomposable initialized aCM bundles of higher rank on $X$ which are stable Ulrich bundles, and studied the properties of them. In the case where $d=4$, E.Coskun and R.S. Kulkarni have constructed a 14-dimensional family of simple Ulrich bundles of rank 2 on a Pfaffian quartic surface $X$ with $c_1=\mathcal{O}_X(3H)$ and $c_2=14$ ([6, Theorem 1.1]). However, in general, if $d\geq4$, it is difficult to consider the problem of whether a given aCM bundle on $X$ splits without special assumptions on the Picard lattice of $X$. In this paper, we will focus on indecomposable initialized aCM bundles of rank 2 on a smooth quartic hypersurface in $\mathbb{P}^3$.

Let $X$ be a smooth quartic hypersurface in $\mathbb{P}^3$. Since $X$ is a K3 surface, the properties of line bundles on K3 surfaces are often used to classify initialized aCM bundles of rank 2 on $X$, by the Chern classes of them. Recently, Gianfranco Casnati has classified indecomposable initialized aCM bundles of rank 2 on $X$, in the case where $X$ is general determinantal ([4]). He has investigated the Chern classes and zero loci of sections of them, by using the fact that the Picard lattice of such a quartic hypersurface is generated by the hyperplane class of $X$ and an Ulrich line bundle on $X$. If an aCM bundle $E$ of rank 2 on $X$ is globally generated, then $\det(E)$ is base point free and big, and the zero locus of a general section of $E$ forms a base point free pencil on a smooth curve $C\in|\det(E)|$. A vector bundle on $X$ becomes globally generated, after twisting it sufficiently many times by the hyperplane class of $X$. Hence, it is interesting to classify aCM bundles of rank 2 on $X$ which are generated by global sections of them in the point of view of the Brill-Noether theory of curves on K3 surfaces.

Let $C$ be a smooth curve of genus $g\geq3$ on $X$, and let $Z$ be a base point free divisor on $C$ with $h^0(\mathcal{O}_C(Z))=r\;(r\geq2)$. Then we obtain a vector bundle of rank $r$ on $X$ as the kernel of the evaluation map 
$$ev:H^0(\mathcal{O}_C(Z))\otimes\mathcal{O}_X\rightarrow \mathcal{O}_C(Z).$$
We denote the dual of it by $E_{C,Z}$, and call it the {\it{Lazarsfeld-Mukai}} bundle associated with $C$ and $Z$. By the construction, such a bundle $E_{C,Z}$ is given by the extension of the rank 1 torsion free sheaf $\mathcal{O}_X(C)\otimes\mathcal{J}_{Z|X}$ by a trivial bundle on $X$, where $\mathcal{J}_{Z|X}$ is the ideal sheaf of $Z$ in $X$. It is well known that $E_{C,Z}$ is globally generated off the base points of $|K_C\otimes\mathcal{O}_C(-Z)|$, and satisfies $h^1(E_{C,Z})=h^2(E_{C,Z})=0$ (for example, see [2, Proposition 2.1], [5, Lemma 1.1]). Conversely, if a vector bundle $E$ which is generated by the global sections satisfies $h^1(E)=h^2(E)=0$, it is of type $E_{C,Z}$ ([5, Lemma 1.2]). In particular, if a rank 2 aCM bundle $E$ on $X$ is globally generated, then $E$ satisfies such conditions, and hence, there exist a smooth curve $C\in|\det(E)|$ and a base point free divisor $Z$ on $C$ with $h^0(\mathcal{O}_C(Z))=2$ such that $E=E_{C,Z}$. In this paper, we will give a necessary condition for a rank 2 vector bundle of type $E_{C,Z}$ on $X$ to be an indecomposable initialized aCM bundle. Our main theorem is as follows.

\newtheorem{thm}{Theorem}[section]

\begin{thm} Let $X$ be a smooth quartic hypersurface in $\mathbb{P}^3$, let $h$ be the hyperplane class of $X$, and let $B$ be a non-trivial initialized aCM line bundle on $X$. Let $C$ be a smooth curve on $X$ of genus $g\geq3$ such that $\mathcal{O}_X(C)\in\mathbb{Z}h\oplus\mathbb{Z}B$, and let $Z$ be a base point free divisor on $C$ with $h^0(\mathcal{O}_C(Z))=2$. If $E_{C,Z}$ is indecomposable initialized and aCM, then $\mathcal{O}_X(C)$ is aCM.
\end{thm} 

\noindent If an aCM bundle of rank 2 on $X$ is globally generated, then it is of type $E_{C,Z}$. Hence, Theorem 1.1 gives a necessary and sufficient condition for the first Chern class $c_1(E)$ of an indecomposable initialized aCM bundle $E$ of rank 2 on $X$ to be aCM, under the condition that $E$ is globally generated. We will use the classification of initialized aCM line bundles on $X$ ([15, Theorem 1.1]) to prove Theorem 1.1. Obviously, a quartic hypersurface $X$ as in Theorem 1.1 satisfies $\rk\Pic(X)\geq2$. In particular, if $X$ is determinantal, then $X$ contains a projectively normal smooth curve of genus 3 and degree 6. Since the class of such a curve in $\Pic(X)$ is an Ulrich line bundle on $X$, if $X$ is general in the sense of Noether-Lefschetz, by Theorem 1.1, the first Chern class of an indecomposable initialized aCM bundle of rank 2 on $X$ which is generated by global sections is an aCM line bundle on $X$, and the existence of such an aCM bundle of rank 2 on $X$ is indicated in [4]. 

Our plan of this paper is as follows. In section 2, we recall basic results about vector bundles and line bundles on K3 surfaces. In section 3, we recall our previous work concerning the classification of initialized aCM line bundles on quartic hypersurfaces in $\mathbb{P}^3$. In section 4, we recall several classical facts about Lazarsfeld-Mukai bundles. In section 5, we give the proof of Theorem 1.1. 

$\;$

{\bf{Notations and conventions}}. We work over the complex number field $\mathbb{C}$. In this paper, a curve and a surface are smooth and projective. Let $X$ be a curve or a surface. Then we denote the canonical bundle of $X$ by $K_X$. Let $E$ be a vector bundle on $X$. Then we denote the rank of $E$, the $i$th Chern class of $E$, and the dual of $E$ by $\rk(E)$, $c_i(E)$, and $E^{\vee}$, respectively. For a divisor or a line bundle $L$ on $X$, we denote by $|L|$ the linear system of $L$. 

Let $C$ be a curve. Then the gonality of $C$ is the minimal degree of pencils on $C$. If the genus of $C$ is $g$ and $Z$ is a line bundle or a divisor of degree $d$ on $C$ with $\dim|Z|=r$, we denote the Brill-Noether number associated with $g$, $r$, and $d$ by $\rho(g,r,d)=g-(r+1)(g-d+r)$. 

Let $X$ be a surface. Then we denote the Picard lattice of $X$ by $\Pic(X)$, and call the rank of it the Picard number of $X$. The Hodge index theorem implies that if the Picard number of $X$ is $\rho$, then the signature of $\Pic(X)$ is $(1,\rho-1)$. Note that this means that $D_1^2D_2^2\leq(D_1.D_2)^2$ for two divisors $D_1$ and $D_2$ on $X$ satisfying  $D_1^2>0$ and $D_2^2>0$. Throughout this paper, we denote the arithmetic operator in $\Pic(X)$ by $+$, and denote by $-L$ the dual of a line bundle $L$ on $X$. We call $X$ a K3 surface if $h^1(\mathcal{O}_X)=0$ and $K_X$ is trivial. In particular, a smooth quartic hypersurface in $\mathbb{P}^3$ is a K3 surface. 

Let $X$ be a quartic hypersurface in $\mathbb{P}^3$. For a hyperplane section $H$ of $X$, we denote the class of $H$ in $\Pic(X)$ by $h$. Moreover, for a vector bundle $E$ on $X$, we denote by $E(L)$ the bundle obtained after twisting it by a line bundle $L$ on $X$. In particular, if $L=\mathcal{O}_X(lH)$, we write it simply as $E(l)$. For a subscheme $Z$ of a variety $Y$, we denote the ideal sheaf of $Z$ in $Y$ by $\mathcal{J}_{Z|Y}$.

\section{Vector bundles and line bundles on K3 surfaces}

In this section, we recall some fundamental facts about vector bundles and line bundles on K3 surfaces. Let $X$ be a K3 surface, and let $E$ be a vector bundle on $X$. First of all, the Riemann-Roch theorem for $E$ is described as follows.
$$\chi(E)=2rk(E)+\displaystyle\frac{c_1(E)^2}{2}-c_2(E),$$
where $\chi(E)=h^0(E)-h^1(E)+h^2(E)$. Note that since $K_X\cong\mathcal{O}_X$, the Serre duality for $E$ is given as follows.
$$h^i(E)=h^{2-i}(E^{\vee})\;\;\; (0\leq i\leq 2).$$
If $L$ is a line bundle, then $L^2$ is an even integer. In particular, if $L^2\geq-2$, then $h^0(L)>0$ or $h^2(L)>0$, and moreover, if $L$ is trivial or satisfies $L.H>0$ for an ample divisor $H$ on $X$, then $h^0(L)>0$.

\newtheorem{df}{Definition}[section]

\begin{df} A non-zero effective divisor $D$ on a surface is called {\rm{$m$-connected}} if $D_1.D_2\geq m$, for each effective decomposition $D=D_1+D_2$.\end{df}
\noindent If a non-zero effective divisor $D$ is 1-connected, then $h^0(\mathcal{O}_D)=1$. Hence, by the exact sequence
$$0\longrightarrow \mathcal{O}_X(-D)\longrightarrow \mathcal{O}_X\longrightarrow\mathcal{O}_D\longrightarrow 0,$$
we have $h^1(\mathcal{O}_X(D))=h^1(\mathcal{O}_X(-D))=0$ ([14, Lemma 2.2]). In particular, if $D$ is an irreducible divisor on $X$, then the dimension of $H^0(\mathcal{O}_X(D))$ and the arithmetic genus $p_a(D)$ of $D$ are given as follows. 
$$h^0(\mathcal{O}_X(D))=2+\displaystyle\frac{D^2}{2},\;\;p_a(D)=1+\displaystyle\frac{D^2}{2}.$$ 
If $D$ is a smooth curve, then $p_a(D)$ coincides with the geometric genus of $D$. In particular, any smooth rational curve $D$ on $X$ satisfies $D^2=-2$. 

\newtheorem{prop}{Proposition}[section]

\begin{prop}{\rm{([14], [13, Proposition 3.4])}}. Let $L$ be a line bundle on a K3 surface $X$ such that $h^0(L)>0$. If $|L|$ is base point free and $L^2>0$, then any member of $|L|$ is 2-connected.\end{prop}

The classification of base point free line bundles on K3 surfaces is well known. We recall the following very useful result.

\begin{prop}{\rm{([14, Proposition 2.6], [8, Chapter 2, Proposition 3.10])}} Let $L$ be a line bundle on a K3 surface $X$ such that $|L|\neq\emptyset$. Assume that $|L|$ has no fixed component. Then one of the following cases occurs.

\smallskip

\smallskip

{\rm{(i)}} $L^2>0$ and the general member of $|L|$ is a smooth irreducible curve of genus $\frac{1}{2}L^2+1$.

{\rm{(ii)}} $L^2=0$ and $L\cong\mathcal{O}_X(kF)$, where $k\geq1$ is an integer and $F$ is a smooth curve of genus one. In this case, $h^1(L)=k-1$. \end{prop} 

\noindent Obviously, by Proposition 2.2, a base point free line bundle is numerical effective. If $C$ is an irreducible curve on $X$ with $C^2\geq0$, then $\mathcal{O}_X(C)$ is base point free ([14]). Therefore, by Proposition 2.2, a line bundle $L$ with $|L|\neq\emptyset$ has no base points outside its fixed components. It is well known that an ample line bundle on $X$ is very ample or hyperelliptic (see [8, Chapter 2, Remark 3.8], [14]). By the characterization of hyperelliptic linear systems on $X$, we have the following assertion.

\begin{prop} {\rm{(cf. [10, Theorem 1.1], [11], [14, Theorem 5.2])}} Let $L$ be a numerically effective line bundle with $L^2\geq4$ on a K3 surface $X$. Then $L$ is very ample if and only if the following conditions are satisfied.

\smallskip

\smallskip

{\rm{(i)}} There is no irreducible curve $E$ such that $E^2=0$ and $E.L=1$ or 2.

{\rm{(ii)}} There is no irreducible curve $E$ such that $E^2=2$ and $L\cong\mathcal{O}_X(2E)$.

{\rm{(iii)}} There is no irreducible curve $E$ such that $E^2=-2$ and $E.L=0$. \end{prop}

\noindent We often investigate the movable part of the linear system $|D|$ of an effective divisor $D$ on $X$ by computing the degree $H.D$ of $D$ for a given very ample divisor $H$ on $X$. We recall the following result at the end of this section.

\begin{prop}{\rm{([16, Corollary 2.1])}} {\it{Let $X$ be a K3 surface, and $L$ be a very ample line bundle on $X$. Then the following statements hold.

\smallskip

\smallskip

{\rm{(i)}} Let $D$ be a nonzero divisor on X with $D^2\geq0$ and $L.D>0$. Then $|D|\neq\emptyset$ and $L.D\geq 3$. In particular, if $D^2=0$ and $L.D=3$, then $|D|$ is an elliptic pencil.

\smallskip

\smallskip

{\rm{(ii)}} There is no effective divisor $D$ on $X$ with $D^2=2$ and $L=\mathcal{O}_X(2D)$.}}\end{prop}

\section{ACM bundles on quartic hypersurfaces in $\mathbb{P}^3$}

Let $X$ be a smooth quartic hypersurface in $\mathbb{P}^3$, let $H$ be a hyperplane section of $X$, and let $h$ be the class of it in $\Pic(X)$. In this section, we recall the result about aCM line bundles on $X$ which is obtained in [15], and give some remarks on aCM bundles on $X$.

\begin{prop} {\rm{([15, Theorem 1.1])}} Let $B$ be a non-trivial line bundle on $X$ with $|B|\neq\emptyset$. Then the following conditions are equivalent.

\smallskip

\smallskip

\noindent {\rm{(i)}} $B$ is aCM and initialized.

\noindent {\rm{(ii)}} One of the following cases occurs.

\smallskip

\smallskip

{\rm{(a)}} $B^2=-2$ and $1\leq h.B\leq 3$.

{\rm{(b)}} $B^2=0$ and $3\leq h.B\leq 4$.

{\rm{(c)}} $B^2=2$ and $h.B=5$.

{\rm{(d)}} $B^2=4,\;h.B=6$ and $|B-h|=|2h-B|=\emptyset.$\end{prop}

\noindent {\bf{Remark 3.1}}. Let $B$ be an aCM line bundle on $X$. Then obviously, $-B$ is also aCM. Assume that $B$ is initialized. By Proposition 3.1, if $B^2=-2$ and $1\leq B.h\leq 2$, then $h-B$ is also initialized. If $B^2=2$, $B^2=0$ and $B.h=4$, or $B^2=-2$ and $B.h=3$ then the same is true for $2h-B$. Moreover, if $B^2=4$, then $3h-B$ is also aCM and initialized. By Proposition 2.4 (i), if $B^2=0$ and $B.h=3$, $|B|$ is an elliptic pencil. By the proof of Proposition 3.1, an initialized aCM line bundle $B$ on $X$ with $B^2\geq2$ is base point free, and hence, by Proposition 2.2 (i), a general member of $|B|$ is smooth. 

\smallskip

\smallskip

If a curve $C\subset\mathbb{P}^3$ which is not necessarily smooth satisfies $h^1(\mathcal{J}_{C|\mathbb{P}^3}(tH))=0$ for each $t\in\mathbb{Z}$, $C$ is called an aCM curve. A curve $C\subset\mathbb{P}^3$ is aCM and smooth if and only if $C$ is projectively normal. Moreover, if $C$ is a divisor on $X$, by applying $\otimes\mathcal{O}_{\mathbb{P}^3}(t)$ to the exact sequence
$$0\longrightarrow\mathcal{O}_{\mathbb{P}^3}(-4)\longrightarrow\mathcal{J}_{C|\mathbb{P}^3}\longrightarrow\mathcal{O}_X(-C)\longrightarrow0,\leqno (3.1)$$
it follows that $\mathcal{O}_X(C)$ is aCM if and only if $C$ is an aCM curve in $\mathbb{P}^3$. 

For a vector bundle $E$ on $X$, we set $H_{\ast}^0(E):=\bigoplus_{l\in\mathbb{Z}}H^0(E(l))$. Then $H_{\ast}^0(E)$ is a graded module over the homogeneous coordinate ring of $X$. Moreover, the following assertion is well known.

\begin{prop}{\rm{([3, Theorem 3.1])}}. Let $E$ be an aCM bundle on $X$, and let $\mu(E)$ be the minimal number of generators of $H_{\ast}^0(E)$. Then we get 
$$\mu(E)\leq4\rk(E).$$\end{prop}
\noindent Proposition 3.2 means that if $E$ is an initialized aCM bundle, then $h^0(E)\leq4\rk(E)$. Moreover, an initialized aCM bundle $E$ satisfying $\mu(E)=4\rk(E)$ is called an {\it{Ulrich}} bundle. In particular, if an initialized aCM line bundle $B$ satisfies the condition as in Proposition 3.1 (ii) (d), $B$ is Ulrich. If there exists such a line bundle $B$ on $X$, then by taking the minimal free resolution of $H_{\ast}^0(B)$ as a module over the homogeneous coordinate ring of $\mathbb{P}^3$, it follows that $X$ is determinantal.


\newtheorem{lem}{Lemma}[section]





\section{Rank 2 Lazarsfeld-Mukai bundles on a quartic hypersurface in $\mathbb{P}^3$}

In this section, we recall the fundamental properties of Lazarsfeld-Mukai bundles of rank 2 on K3 surfaces and several important results. Moreover, we remark some facts about rank 2 Lazarsfeld-Mukai bundles on quartic hypersurfaces in $\mathbb{P}^3$. Let $X$, $H$, and $h$ be as in section 3, let $C$ be a curve of genus $g\geq3$ on $X$, and let $Z$ be a base point free divisor of degree $d$ on $C$ with $\dim|Z|=r$. Then by taking the dual of the exact sequence associated with the evaluation map of $\mathcal{O}_C(Z)$, the Lazarsfeld-Mukai bundle $E_{C,Z}$ associated with $C$ and $Z$ fits the following exact sequence.
$$0\longrightarrow H^0(\mathcal{O}_C(Z))^{\vee}\otimes\mathcal{O}_X\longrightarrow E_{C,Z}\longrightarrow K_C\otimes\mathcal{O}_C(-Z)\longrightarrow 0.\leqno(4.1)$$
Moreover, $E_{C,Z}$ has the following properties.

\begin{prop}{\rm{(See [5, Lemma 1.1], [2, Proposition 2.1])}} Let the notations be as above. Then

$\;$

{\rm{(a)}} $c_1(E_{C,Z})=\mathcal{O}_X(C)$.

\smallskip

\smallskip

{\rm{(b)}} $c_2(E_{C,Z})=d$.

\smallskip

\smallskip

{\rm{(c)}} $h^1(E_{C,Z})=h^2(E_{C,Z})=0$.

\smallskip

\smallskip

{\rm{(d)}} $\chi(E_{C,Z}^{\vee}\otimes E_{C,Z})=2(1-\rho(g,r,d))$.

\smallskip

\smallskip

{\rm{(e)}} $h^0(E_{C,Z})=g-d+1+2r$.

\smallskip

\smallskip

{\rm{(f)}} $E_{C,Z}$ is globally generated off the base points of $|K_C\otimes\mathcal{O}_C(-Z)|$. \end{prop}

\noindent In Proposition 4.1, the property (e) follows by easy computation from the exact sequence (4.1). By the property (d), if $\rho(g,r,d)<0$, then $E_{C,Z}$ is not simple. Moreover, we note the following assertion.

\begin{prop} {\rm{([5, Lemma 1.2])}}. Let $E$ be a vector bundle on $X$ generated by its global sections such that $h^1(E)=h^2(E)=0$. Then $E$ is of type $E_{C,Z}$.\end{prop}

From now on we assume that the rank of $E_{C,Z}$ is 2. Then the following result is well known. 

\begin{prop} {\rm{([7, Lemma 4.4])}}. Let the notations be as above. If $E_{C,Z}$ is not simple, then there exist two line bundles $M$ and $N$ on $X$ and a 0-dimensional subscheme $Z^{'}\subset X$ of finite length such that

$\;$

{\rm{(a)}} $h^0(M)\geq2,\;h^0(N)\geq2;$

\smallskip

\smallskip

{\rm{(b)}} $N$ is base point free;

\smallskip

\smallskip

{\rm{(c)}} There exists an exact sequence
$$0\longrightarrow M\longrightarrow E_{C,Z}\longrightarrow N\otimes\mathcal{J}_{Z^{'}|X}\longrightarrow0.$$
Moreover, if $h^0(M-N)=0$, then the length of $Z^{'}$ is zero. \end{prop}

\noindent {\bf{Remark 4.1}}. If $E_{C,Z}$ is not simple, the exact sequence as in Proposition 4.3 is obtained by taking the image and the kernel of an endomorphism of $E_{C,Z}$ dropping the rank everywhere. Moreover, in Proposition 4.3, if $h^0(M-N)>0$, then $C.(M-N)\geq0$, and if $h^0(M-N)=0$, then $E_{C,Z}\cong N\oplus M$. Hence, we can assume that $C.M\geq C.N$, that is, $M^2\geq N^2$.

\smallskip

\smallskip

\noindent In particular, in Proposition 4.3, if $|Z|$ is a gonality pencil, then we have the following more detailed assertion.

\begin{prop}{\rm{([5, Proposition 2.3])}}. Let the notations be as above. Assume that $d$ is the minimal gonality of a smooth curve in $|\mathcal{O}_X(C)|$, and $\mathcal{O}_C(Z)$ is a line bundle corresponding to a $g_d^1$ on $C$. Then one of the following cases occurs.

\smallskip

\smallskip

\noindent{\rm{(a)}} $E_{C,Z}$ is simple and then $\rho(g,1,d)\geq0$.

\smallskip

\smallskip

\noindent{\rm{(b)}} $\rho(g,1,d)=0$, and there exist a base point free line bundle $N$ on $X$, a $(-2)$-curve $\Delta$ on $X$ with $N.\Delta=1$, and a point $P$ of $X$ such that $C\in |N^{\otimes2}(\Delta)|$, and there exists the following exact sequence.
$$0\longrightarrow N(\Delta)\longrightarrow E_{C,Z}\longrightarrow N\otimes\mathcal{J}_{P|X}\longrightarrow0.$$

\smallskip

\smallskip

\noindent{\rm{(c)}} There exist two line bundles $M$ and $N$ on $X$ satisfying the following properties.

\smallskip

{\rm{($\text{c}_1$)}} $h^0(M)\geq2$ and $h^0(N)\geq2$.

\smallskip

{\rm{($\text{c}_2$)}} $N$ is base point free.

\smallskip

{\rm{($\text{c}_3$)}} $h^1(M)=h^1(N)=0$.

\smallskip

{\rm{($\text{c}_4$)}} The {\rm{(}}possibly empty{\rm{)}} fixed component $\Delta$ of $M$ satisfies $\Delta.C\leq1$. Moreover, if $\Delta.C=1$, then the gonality of the smooth curves in $|\mathcal{O}_X(C)|$ is constant. If $C$ is ample, then $\Delta$ is empty.

\smallskip

{\rm{($\text{c}_5$)}} There exists the following exact sequence.
$$0\longrightarrow M\longrightarrow E_{C,Z}\longrightarrow N\longrightarrow0.$$

\end{prop}

\noindent In Proposition 4.4, if $E_{C,Z}$ is not simple, one of the cases (b) and (c) occurs. By the proof of Proposition 4.4, if $\rho(g,1,d)<0$, the case (c) occurs. In this case, we may also assume that $M^2\geq N^2$, by Remark 4.1.

$\;$

\noindent{\bf{Remark 4.2}}. Since $h^0(\mathcal{O}_C(Z))=2$, by applying $\otimes\mathcal{O}_X(-l)\;(l\in\mathbb{Z})$ to the exact sequence (4.1), we have
$$\chi(E_{C,Z}(-l))=4l^2-lC.H+g+3-d.$$

\smallskip

\smallskip

\noindent {\bf{Remark 4.3}}. Assume that $E_{C,Z}$ is initialized and aCM. Then for $l\geq1$,
$$\chi(E_{C,Z}(-l))=h^2(E_{C,Z}(-l))=h^0(E_{C,Z}(lH-C)).\leqno(4.3)$$
By applying $\otimes\mathcal{O}_X(lH-C)$ to the exact sequence (4.1), we have the following exact sequence.
$$0\longrightarrow H^0(\mathcal{O}_C(Z))^{\vee}\otimes\mathcal{O}_X(lH-C)\longrightarrow E_{C,Z}(lH-C)\longrightarrow \mathcal{O}_C(lH)\otimes\mathcal{O}_C(-Z)\longrightarrow 0.\leqno(4.2)$$
If $C$ is an aCM curve in $\mathbb{P}^3$, then $h^1(\mathcal{O}_X(lH-C))=0$, and hence, by (4.2) and the exact sequence
$$0\longrightarrow \mathcal{O}_X(lH-C)\longrightarrow\mathcal{O}_X(lH)\otimes\mathcal{J}_{Z|X}\longrightarrow\mathcal{O}_C(lH)\otimes\mathcal{O}_C(-Z)\longrightarrow 0,$$
we have $h^0(\mathcal{O}_X(lH)\otimes\mathcal{J}_{Z|X})=\chi(E_{C,Z}(-l))-h^0(\mathcal{O}_X(lH-C))$.

\smallskip

\smallskip

\noindent Remark 4.3 means that if $E_{C,Z}$ is initialized and aCM and $C$ is an aCM curve in $\mathbb{P}^3$, then we can determine the Hilbert function of the ideal of $Z$ in $X$.

\smallskip

\smallskip

\noindent {\bf{Remark {4.4}}}. By the exact sequence
$$0\longrightarrow\mathcal{O}_X(-1)\longrightarrow\mathcal{O}_X(C-H)\otimes\mathcal{J}_{Z|X}\longrightarrow\mathcal{O}_C(C-H)\otimes\mathcal{O}_C(-Z)\longrightarrow0,$$
applying $\otimes\mathcal{O}_X(-1)$ to (4.1), we have $h^0(E_{C,Z}(-1))=h^0(\mathcal{O}_X(C-H)\otimes\mathcal{J}_{Z|X})$. Therefore, $E_{C,Z}$ is initialized if and only if $h^0(\mathcal{O}_X(C-H)\otimes\mathcal{J}_{Z|X})=0$.

\smallskip

\smallskip

At the end of this section, we prepare the following lemma.

\begin{lem} Let the notations be as above. If $E_{C,Z}$ is initialized and aCM, then $C.H\leq12$.\end{lem}

\noindent {\it{Proof}}. By Proposition 4.1 (e), we have $h^0(E_{C,Z})=g+3-d$. Hence, by Proposition 3.2, we have $g-5\leq d$. Since $h^0(E_{C,Z}(-1))=h^1(E_{C,Z}(-1))=0$, by Remark 4.2, we have $h^2(E_{C,Z}(-1))=g+7-d-C.H.$ Since
$$d=g+7-C.H-h^2(E_{C,Z}(-1))\leq g+7-C.H,$$
we have $C.H\leq12$. $\hfill\square$

\section{Proof of Theorem 1.1}

In this section, we give the proof of Theorem 1.1. Let $X\subset\mathbb{P}^3$ be a smooth quartic hypersurface. Let $B$ be an initialized aCM line bundle on $X$ which is not trivial, let $C$ be a curve of genus $g\geq3$ on $X$ such that $\mathcal{O}_X(C)\in\mathbb{Z}h\oplus\mathbb{Z}B\subset\Pic(X)$, and let $Z$ be a base point free divisor of degree $d$ on $C$ such that $h^0(\mathcal{O}_C(Z))=2$. Assume that $E_{C,Z}$ is indecomposable initialized and aCM. First of all, we prepare the following lemma.

\begin{lem} Let $s$ and $t$ be integers satisfying $C\in |sh+tB|$. Assume that $|t|\geq2$. Then we obtain the following assertions.

$\;$

\noindent {\rm{(i)}} Assume that $B^2=-2$.

\smallskip

{\rm{(a)}} If $B.h=1$, then $(s,t)=(3,-2)$.

\smallskip

{\rm{(b)}} If $B.h=2$, then $(s,t)=(2,2)$ or $(4,-2)$.

\smallskip

{\rm{(c)}} If $B.h=3$, then $(s,t)=(4,-2)$.

\smallskip

\smallskip

\noindent {\rm{(ii)}} If $B^2=0$ and $B.h=4$, then $(s,t)=(1,2)$ or $(5,-2)$.

\smallskip

\smallskip

\noindent {\rm{(iii)}} If $B^2=4$, then $(s,t)=(0,2)$ or $(6,-2)$. \end{lem}

\noindent {\it{Proof}}. Assume that $B^2=-2$. We consider the case where $B.h=1$. By Lemma 4.1, we have $C.H=4s+t\leq12$. Since $C.B=s-2t\geq0$, we have $9t\leq12$, and hence, we have $t\leq1$. On the other hand, since $(h-B)^2=0$ and $h.(h-B)=3$, by Remark 3.1, $|h-B|$ is an elliptic pencil. Since $C^2>0$, by Proposition 2.1, we have $C.(h-B)=3(s+t)>0$, and hence, $s\geq 1-t$. Since $4s+t\leq12$, we have $t\geq-2$. By the assumption, $(s,t)=(3,-2)$.

We consider the case where $B.h=2$. By the same reason as above, we have $C.H=4s+2t\leq12$. Since $C.B=2s-2t\geq0$, we have $t\leq2$. On the other hand, since $(h-B)^2=-2$, $h.(h-B)=2$, and $C\in |(s+t)h-t(h-B)|$, $h-B$ plays the same role of $B$. By the same reason as above, we have $-t\leq2$. Hence, by the assumption, we have $(s,t)=(2,2)$ or $(4,-2)$.

We consider the case where $B.h=3$. Since 
$$C.H=4s+3t\leq12\text{ and }C.B=3s-2t\geq0,\leqno(5.1)$$
we have $t\leq2$. On the other hand, since $(2h-B)^2=2$ and $h.(2h-B)=5$, $|2h-B|\neq\emptyset$. Since $g\geq3$, by the Hodge index theorem, we have 
$$C.(2h-B)=5s+8t\geq3.$$
Hence, we have $-17t\leq 48$, and hence, $t\geq-2$. If $t=2$, there is no integer $s$ satisfying the inequality (5.1). By the assumption, we have $(s,t)=(4,-2)$.

Assume that $B^2=0$ and $B.h=4$. Since $C^2>0$ and the movable part of $|B|$ is not empty, by Proposition 2.1, we have $C.B=4s>0$, and hence, $s\geq1$. Since $C.H=4s+4t\leq12$, we have $t\leq2$. Since $(2h-B)^2=0$, $h.(2h-B)=4$, and $C\in|(s+2t)h-t(2h-B)|$, $2h-B$ plays the same role of $B$. By the same reason as above, we have $t\geq-2$. By the assumption, we have $(s,t)=(1,2)$ or $(5,-2)$.

Assume that $B^2=4$. By Proposition 3.1, $B.h=6$. Hence, by Lemma 4.1, we have $C.H=4s+6t\leq12$. Since $g\geq3$, by the Hodge index theorem, we have $C.B=6s+4t\geq4$, and hence, $t\leq2$. On the other hand, since $(3h-B)^2=4$, $h.(3h-B)=6$, and $C\in|(s+3t)h-t(3h-B)|$, by the same reason as above, we have $t\geq-2$. By the assumption, we have $(s,t)=(0,2)$ or $(6,-2)$. $\hfill\square$

\smallskip

\smallskip

Before the proof of Theorem 1.1, we recall the following definition.

\begin{df} Let $L$ be a line bundle on $X$. We say that a locally complete intersection subscheme $W\subset X$ of dimension 0 is {\rm{Cayley-Bacharach}} {\rm{(CB}} for short{\rm{)}} with respect to $L$ if, for each $W^{'}\subset W$ of degree $\deg(W)-1$, the natural morphism $H^0(\mathcal{J}_{W|X}\otimes L)\rightarrow H^0(\mathcal{J}_{W^{'}|X}\otimes L)$ is an isomorphism. \end{df}

\noindent The following result on the construction of vector bundles of rank 2 on $X$ is well known.

\begin{prop}{\rm{([9, Theorem 5.1.1])}}. Let $L$ be a line bundle on $X$ and let $W\subset X$ be a locally complete intersection subscheme of dimension 0. Then there exists a vector bundle $E$ of rank 2 on $X$ which fits the exact sequence
$$0\longrightarrow\mathcal{O}_X\longrightarrow E\longrightarrow L\otimes\mathcal{J}_{W|X}\longrightarrow 0$$
if and only if $W$ is CB with respect to $L$.\end{prop}

\noindent Since $Z$ is base point free, it is CB with respect to $\mathcal{O}_X(C)$. Hence, by Proposition 5.1, we can construct a vector bundle $E$ of rank 2 on $X$ with $c_1(E)=\mathcal{O}_X(C)$ which has a global section whose zero locus is $Z$. In particular, by taking a global section of $E_{C,Z}$ which has the same zero locus $Z$, $E_{C,Z}$ fits the exact sequence
$$0\longrightarrow\mathcal{O}_X\longrightarrow E_{C,Z}\longrightarrow\mathcal{O}_X(C)\otimes\mathcal{J}_{Z|X}\longrightarrow0.\leqno(5.2)$$

\smallskip

\smallskip

\noindent{\it{Proof of Theorem 1.1}}. If $B^2=0$ and $B.h=3$, then we have $(h-B)^2=-2$ and $(h-B).h=1$. Moreover, if $B^2=2$ and $B.h=5$, then we have $(2h-B)^2=-2$ and $(2h-B).h=3$. Hence, by Proposition 3.1, it is sufficient to show that the each case as in Lemma 5.1 does not occur.

We consider the case where $B^2=-2$. Assume that $B.h=1$ and $C\in|3h-2B|$. Let $F=h-B$. Since $h.F=3$, $F^2=0$, and $h$ is very ample, by Proposition 2.4 (i), $|F|$ is an elliptic pencil. Since $E_{C,Z}$ is initialized, we have $h^0(E_{C,Z}(-h-F))=0$. Since $B.h=1$, the member of $|B|$ is a $(-2)$-curve. Let $\Gamma\in |B|$. Then we have the following exact sequence.
$$0\longrightarrow E_{C,Z}(-2)\longrightarrow E_{C,Z}(-h-F)\longrightarrow E_{C,Z}(-h-F)|_{\Gamma}\longrightarrow0.$$
Since $h^0(E_{C,Z}(-2))=h^1(E_{C,Z}(-2))=0$, we have $h^0(E_{C,Z}(-h-F)|_{\Gamma})=0$. Since $\deg(E_{C,Z}(-h-F)|_{\Gamma})=-H.\Gamma=-1$, there exists an integer $a$ such that $E_{C,Z}(-h-F)|_{\Gamma}\cong\mathcal{O}_{\Gamma}(-1+a)\oplus\mathcal{O}_{\Gamma}(-a)$. Hence, we have $h^0(E_{C,Z}(-h-F)|_{\Gamma})>0$. This is a contradiction.

Assume that $B.h=2$ and $C\in |2h+2B|$. Since $g=13$, by Proposition 4.1 (e), we have $h^0(E_{C,Z})=16-d$. By Proposition 3.2, we have $d\geq8$. Since $h^1(E_{C,Z}(-1))=0$ and $C.H=12$, by Remark 4.2, we have $\chi(E_{C,Z}(-1))=8-d\geq0$. Hence, we have $d=8$. Since $\mathcal{O}_X(C-H)=2B+h$ and, by Remark 4.4, $h^0(\mathcal{O}_X(C-H)\otimes\mathcal{J}_{Z|X})=0$, we have $h^0((B+h)\otimes\mathcal{J}_{Z|X})=0$. By applying $\otimes \mathcal{O}_X(-h-B)$ to the exact sequence (5.2), we have the following exact sequence
$$0\longrightarrow\mathcal{O}_X(-h-B)\longrightarrow E_{C,Z}(-h-B)\longrightarrow \mathcal{O}_X(h+B)\otimes\mathcal{J}_{Z|X}\longrightarrow0.$$
Since $h^1(-h-B)=h^0(-h-B)=0$, we have 
$$h^2(E_{C,Z}(-h-B))=h^0(E_{C,Z}(-h-B))=0.$$
Since $\det(E_{C,Z}(-h-B))=\mathcal{O}_X$ and $c_2(E_{C,Z}(-h-B))=2$. Hence, we obtain the contradiction $h^1(E_{C,Z}(-h-B))=-\chi(E_{C,Z}(-h-B))=-2$.

Assume that $C\in|4h-2B|$. Since $(h-B)^2=-2$, $h.(h-B)=2$, and $4h-2B=2h+2(h-B)$, $h-B$ plays the same role of $B$. Hence, by the same reason as above, this case does not occur.

Assume that $B.h=3$ and $C\in|4h-2B|$. Since $h^1(E_{C,Z}(-1))=0$, $g=5$, and $C.H=10$, by Remark 4.2, $\chi(E_{C,Z}(-1))=2-d\geq0$. Therefore, we have $d=2$. Obviously, $d$ is the minimal gonality of smooth curves in $|C|$ and $\rho(g,1,d)=-3<0$. Hence, by Proposition 4.4 and Remark 4.1, there exist line bundles $M,N\in\Pic(X)$ such that $h^0(M)\geq2$, $h^0(N)\geq2$, $h^1(M)=h^1(N)=0$, $M^2\geq N^2$, $|N|$ is base point free, and $E_{C,Z}$ fits the exact sequence
$$0\longrightarrow M\longrightarrow E_{C,Z}\longrightarrow N\longrightarrow 0.$$

Assume that $N^2=0$. Then $|N|$ is an elliptic pencil. Since $(2h-B).h=5$ and $(2h-B)^2=2$, by Remark 3.1, $|2h-B|$ is base point free. Hence, by Proposition 2.1, we have $(2h-B).N\geq2$. Therefore, we have $M.N=C.N\geq4$. However, this contradicts the facts that $d=2$.

Assume that $N^2=2$. Since $2=M.N=C.N-2$, we have $(2h-B).N=2$, and hence, $(2h-B-N)^2=0$. If $2h-B\neq N$, the movable part of $|2h-B-N|$ or $|N-2h+B|$ is not empty. Since $N.(2h-B-N)=0$, this contradicts the assertion of Proposition 2.1. Hence, $N=2h-B$. This means that $E_{C,Z}\cong\mathcal{O}_X(2h-B)^{\oplus 2}$. However, this contradicts the assumption that $E_{C,Z}$ is indecomposable.

Since $M^2\geq N^2$, by the Hodge index theorem, if $N^2\geq4$, we have the contradiction $2=M.N\geq N^2\geq4$. Therefore, this case does not occur.

We consider the case where $B^2=0$. Assume that $B.h=4$ and $C\in|h+2B|$. Since $g=11$, by the same reason as above, we have $h^0(E_{C,Z})=14-d\leq8$, and hence, $d\geq6$. Since $C.H=12$ and $h^1(E_{C,Z}(-1))=0$, by Remark 4.2, we have $\chi(E_{C,Z}(-1))=6-d\geq0$, and hence, $d=6$. Since $\rho(g,1,d)=-1<0$, by Proposition 4.3 and Remark 4.1, there exist line bundles $M,N\in\Pic(X)$, and 0-dimensional subscheme $Z^{'}\subset X$ such that $h^0(M)\geq2$, $h^0(N)\geq2$, $M^2\geq N^2$, $|N|$ is base point free, and $E_{C,Z}$ fits the exact sequence
$$0\longrightarrow M\longrightarrow E_{C,Z}\longrightarrow N\otimes\mathcal{J}_{Z^{'}|X}\longrightarrow0.$$

Assume that $N^2=0$. Since $|N|$ is base point free, by Proposition 2.2 (ii), there exist a positive integer $r$ and an elliptic curve $F$ on $X$ such that $|N|=|rF|$. By Proposition 2.3, $H.F\geq3$, and hence, 
$$3r\leq r(h+2B).F=C.N=M.N\leq6.$$
If $r\geq2$, then we have $r=2$, and hence, the length of $Z^{'}$ is zero. However, since $h^1(E_{C,Z})=0$ and $h^2(M)=h^0(-M)=0$, we have $h^1(N)=0$. This is a contradiction. Hence, $r=1$. Assume that $N.B\leq1$. Then we have $h^0(\mathcal{O}_F(B))\leq1$. If $h^0(B-N)=0$, then by the exact sequence
$$0\longrightarrow\mathcal{O}_X(B-N)\longrightarrow\mathcal{O}_X(B)\longrightarrow\mathcal{O}_F(B)\longrightarrow 0,$$
we have $h^0(\mathcal{O}_X(B))\leq1$. However, since $B^2=0$, the movable part of $|B|$ is not empty. This is a contradiction. Since $h^0(B-N)>0$, $|M(-1)|\neq\emptyset$. This means that $h^0(E_{C,Z}(-1))>0$. This contradicts the assumption that $E_{C,Z}$ is initialized. Therefore, $N.B\geq2$. Since $N.h\geq3$, we have 
$$7\leq N.(h+2B)=N.C=N.M.$$
This contradicts the fact that $d=6$.

Assume that $N^2>0$. Since $M^2\geq N^2$, by the Hodge index theorem, we have $M.N\geq N^2$. Hence, $4N^2\leq (M+N)^2=C^2=20$. Since $N^2$ is even, we have $N^2\leq4$. Assume that $N^2=2$. Since the movable part of $|B|$ is not empty, by Proposition 2.1, we have $N.B\geq2$. By the Hodge index theorem, we have $N.h\geq3$.  Hence, we have
$$5\leq N.(h+2B)-2=N.C-2=N.M\leq6.$$
Hence, we have $N.B=2$ and $N.h=3$ or 4. We note that $(N-B)^2=-2$. If $N.h=4$, we have $h.(N-B)=0$. However, this contradicts the ampleness of $h$. Assume that $N.h=3$. Since $h.(B-N)=1$, the member of $|B-N|$ is a $(-2)$-curve. If we let $\Gamma\in|B-N|$, we have $B=\mathcal{O}_X(N+\Gamma)$ and $N.\Gamma=N.(B-N)=0$. This contradicts the fact that $h^1(B)=0$. 

Assume that $N^2=4$. Then we have $C.N-4=M.N\leq6$, and hence, $C.N\leq10$. By the Hodge index theorem, we have $C.N\geq\sqrt{80}$, and hence, we have $C.N=10$ or 9. Assume that $C.N=10$. Since, by the Hodge index theorem, $h.N\geq4$, we have $2B.N=10-h.N\leq6$. Since, by the same reason as above, $B.N\geq2$, we have $(B.N,h.N)=(2,6)$ or $(3,4)$. If $B.N=2$, we have
$$(N-B)^2=4-2N.B=0\text{ and }h.(N-B)=2.$$ 
This contradicts the assertion of Proposition 2.4 (i). If $B.N=3$, we have
$$(N-B)^2=-2\text{ and }h.(N-B)=0.$$
This contradicts the ampleness of $h$. Assume that $C.N=9$. By the same reason as above, we have $2B.N=9-h.N\leq5$. Hence, we have $B.N=2$ and $h.N=5$. Then we have $(N-B)^2=0$ and $h.(N-B)=1$. However, by the same reason as above, this is a contradiction. By the above argument, the case that $B.h=4$ and $C\in|h+2B|$ does not occur. 

Assume that $B.h=4$ and $C\in|5h-2B|$. Since $(2h-B)^2=0$, $h.(2h-B)=4$, and $5h-2B=h+2(2h-B)$, $2h-B$ plays the same role of $B$. Hence, this case also does not occur.

Here, we prepare the following lemma.

\begin{lem} If $B^2=4$, then the minimal gonality of curves in $|2B|$ is $4$.  \end{lem}

\noindent {\it{Proof of Lemma 5.2}}. Let $d_0$ be the minimal gonality of curves in $|2B|$. Let $C_0\in |2B|$ be a smooth curve of gonality $d_0$, and let $Z_0$ be a divisor of degree $d_0$ on $C_0$ such that $|Z_0|$ is a pencil on $C_0$. Note that the genus of $C_0$ is 9. First of all, $h^1(B)=h^0(-B)=0$ and, by Remark 3.1, $|B|$ is base point free. By Proposition 4.2, there exist a smooth curve $C_1\in|2B|$ and a base point free pencil $Z_1$ of degree 4 on $C_1$ such that $E_{C_1,Z_1}=B^{\oplus 2}$. Hence, $d_0\leq4$. Assume that $d_0\leq3$. Since $\rho(9,1,d_0)\leq-5<0$, by Proposition 4.4 and Remark 4.1, there exist line bundles $M,N\in\Pic(X)$ such that $h^0(M)\geq2$, $h^0(N)\geq2$, $h^1(N)=h^1(M)=0$, $M^2\geq N^2$, $|N|$ is base point free, and $E_{C_0,Z_0}$ fits the exact sequence
$$0\longrightarrow M\longrightarrow E_{C_0,Z_0}\longrightarrow N\longrightarrow 0.$$

Assume that $N^2=0$. Since, by Remark 3.1, $B$ is base point free, we have $B.N\geq2$ by Proposition 2.1. Hence, $M.N=C.N=2B.N\geq4$. This contradicts the assumption that $d_0\leq3$. 

Assume that $N^2=2$. By the Hodge index theorem, we have $B.N\geq3$. Hence, we have the contradiction $M.N=C.N-2\geq4$.

Since $M^2\geq N^2$, if $N^2\geq4$, by the Hodge index theorem, we have $M.N\geq4$. This is also a contradiction. Hence, we have $d_0=4$.

$\;$

W consider the case where $B^2=4$. Assume that $C\in|2B|$. Since $g=9$, by the same reason as above, we have $h^0(E_{C,Z})=12-d\leq8$, and hence, $d\geq4$. Since $h^1(E_{C,Z}(-1))=0$, and $C.H=12$, we have $\chi(E_{C,Z}(-1))=4-d\geq0$. Therefore, we have $d=4$.

Since $\rho(9,1,4)=-3<0$, by Lemma 5.2, Proposition 4.4, and Remark 4.1, there exist line bundles $M,N\in\Pic(X)$ such that $h^0(M)\geq2$, $h^0(N)\geq2$, $h^1(N)=h^1(M)=0$, $M^2\geq N^2$, $|N|$ is base point free, and $E_{C,Z}$ fits the exact sequence
$$0\longrightarrow M\longrightarrow E_{C,Z}\longrightarrow N\longrightarrow 0.\leqno(5.3)$$
Since by the assumption, $E_{C,Z}$ is indecomposable, $M\neq N$ and by Remark 4.1, $h^0(M-N)>0$. Hence, we have $h.M>h.N$. Since $2h.N<H.C=12$, we have $h.N\leq5$. On the other hand, since $2B.N=C.N$ and $C.N-N^2=M.N=4$, we have $(B-N)^2=0$. Since $h.(B-N)\geq1$, $|B-N|$ is not empty and the movable part of it is not empty. Therefore, we have $N^2\leq2$. In fact, if $N^2\geq4$, by the Hodge index theorem, we have $h.N\geq4$. Hence, we have $h.(B-N)\leq2$. This contradicts the assertion of Proposition 2.4 (i).

Assume that $N^2=0$. Since $h^1(M)=h^1(N)=0$, we have 
$$\displaystyle\frac{M^2}{2}+2=h^0(M)=h^0(E_{C,Z})-h^0(N)=6,$$
and hence, we obtain $M^2=8$. Since by Proposition 2.4 (i), $h.N\geq3$, we have $h.(B-N)\leq3$. Since the movable part of $|B-N|$ is not empty, we have $h.(B-N)=3$, that is, $h.N=3$. Hence, we have $h.M=H.C-3=9$. Since $(M-h)^2=-6$, we have $h^1(M(-1))\neq0$. On the other hand, since $h.(N-h)=-1<0$, $|N(-1)|$ is empty. Since $h^1(E_{C,Z}(-1))=0$, by the exact sequence which is obtained by applying $\otimes\mathcal{O}_X(-1)$ to (5.3) 
$$0\longrightarrow M(-1)\longrightarrow E_{C,Z}(-1)\longrightarrow N(-1)\longrightarrow 0,$$
we have $h^1(M(-1))=0$. However, this is a contradiction.

Assume that $N^2=2$. Since $C.N-2=M.N=4$, we have $B.N=3$. Hence, we have $(B-N).N=1$. However, since by Remark 3.1, $|B|$ is base point free and $|B-N|\neq\emptyset$, this contradicts Proposition 2.1. Hence, the case that $C\in|2B|$ does not occur. Assume that $C\in|6h-2B|$. Since $(3h-B)^2=4$, $(3h-B).h=6$, and $|2h-B|=|B-h|=\emptyset$, $3h-B$ plays the same role of $B$. Hence, by the same argument as above, we have the contradiction. $\hfill\square$

\smallskip

\smallskip

We can find an example of a rank 2 initialized aCM bundle of type $E_{C,Z}$ on a quartic hypersurface in $\mathbb{P}^3$ given by a non-split extension of an aCM line bundle such that $C$ is not aCM . 

\smallskip

\smallskip

\noindent {\bf{Example 5.1}}. Let $g(x_1,x_2,x_3)$ be a general non-degenerate ternary quartic form, and let $X_g$ be the smooth quartic hypersurface in $\mathbb{P}^3$ defined by the equation ${x_0}^4-g(x_1,x_2,x_3)=0$. Let $\pi_1:S\rightarrow\mathbb{P}^2$ be the double covering of $\mathbb{P}^2$ branched along the smooth quartic curve on $\mathbb{P}^2$ defined by $g(x_1,x_2,x_3)=0$. Then $S$ is a DelPezzo surface of degree 2, and $X_g$ is obtained as the double covering $\pi_2:X_g\rightarrow S$ of $S$ branched along a smooth curve $H\in|-2K_S|$. Since $S$ is described as the blow-up $\varphi:S\rightarrow\mathbb{P}^2$ of $\mathbb{P}^2$ at seven points $P_1,\cdots,P_7$ in general position, if we let $l=\varphi^{\ast}\mathcal{O}_{\mathbb{P}^2}(1)$, let $E_i=\varphi^{-1}(P_i)$ and let $e_i$ be the class of $E_i$ in $\Pic(S)$ for $1\leq i\leq7$, then $\Pic(S)$ is generated by $e_i\;(1\leq i\leq7)$ and $l$. Hence, $\Pic(X_g)$ is generated by ${\pi_2}^{\ast}e_i\;(1\leq i\leq7)$ and ${\pi_2}^{\ast}l$ (see [6, Proposition 4.1]). Let $\tilde{l}={\pi_2}^{\ast}l$, let $\tilde{e_i}={\pi_2}^{\ast}e_i\;(1\leq i\leq7)$, and let $h$ be the class of ${\pi_2}^{-1}(H)$ in $\Pic(X_g)$. Then $h$ is the hyperplane class of $X_g$ and $h=3\tilde{l}-\sum_{1\leq i\leq7}\tilde{e_i}$. Let $f=2\tilde{l}-\sum_{1\leq i\leq4}\tilde{e_i}$ and let $f_j=\tilde{l}-\tilde{e_j}\;(1\leq j\leq7)$. Since $f^2={f_j}^2=0$ and $h.f=h.f_j=4$, and $\Pic(X_g)$ is an even lattice, by Proposition 2.3 and Proposition 2.4, $|f|$ and $|f_j|\;(1\leq j\leq7)$ are elliptic pencils. Hence, a vector bundle $E$ of rank 2 on $X_g$ which fits the exact sequence
$$0\longrightarrow\mathcal{O}_X(f)\longrightarrow E\longrightarrow\mathcal{O}_X(f_j)\longrightarrow0\leqno(5.4)$$
is of type $E_{C,Z}$. Assume that $j\in\{5,6,7\}$. Since $(f-f_j)^2=-8$, $h^1(f-f_j)\neq0$. Hence, we can take $E$ such that the exact sequence (5.4) does not split. Then obviously, such a vector bundle $E$ is initialized and aCM, since $f$ and $f_j$ are initialized and aCM by Proposition 3.1. However, since $(f+f_j-2h)^2=-8$, we have $h^1(f+f_j-2h)\neq0$. This means that $c_1(E)$ is not aCM.

$\;$

\noindent {\bf{Acknowledgements}}

\smallskip

\smallskip

\noindent I would like to thank Akira Ohbuchi and Komeda Jiryo for the opportunity to talk about this topic at the 17 th symposium on algebraic curves theory.

\end{document}